\documentclass[11pt]{amsart}
\usepackage{amsmath}
\usepackage{amssymb}
\usepackage{amsthm}

\newcommand{\QQ}{\ensuremath{\mathbb Q\,}}
\newcommand{\ZZ}{\ensuremath{\mathbb Z}}
\newcommand{\FF}{\ensuremath{\mathbb F}}

\newcommand{\NN}{\ensuremath{\mathbb N}}

\newcommand{\eps}{\varepsilon}
\newcommand{\ord}{\mathrm{ord}}

\newcommand{\Gal}{\mathrm{Gal}}
\newcommand{\Frob}{\text{Frob}}
\newcommand{\diag}{\text{diag}}
\renewcommand{\P}{\ensuremath{\mathfrak{P}}}
\newcommand{\p}{\ensuremath{\mathfrak{p}}}

\renewcommand{\o}{\ensuremath{\mathfrak{o}}}

\newcommand{\norm}{\mathbf{N}}

\newcommand{\bF}{\bar\FF}
\newcommand{\bFp}{\bF_p}
\newcommand{\tK}{\tilde K}
\newcommand{\tG}{\tilde G}
\newcommand{\tk}{\tilde k}

\newtheorem{lemma}{Lemma}

\newtheorem{prop}{Proposition}

\title{$A_6$-extensions of \QQ\ and the mod $p$ Cohomology of $GL_3(\ZZ)$} 
\author{A.\ Ash${}^1$}
\author{D.\ Pollack }
\author{ W.\ Sinnott}
\date{November 30, 2003; revised October 19, 2004}
\begin{document}
\maketitle\footnotetext[1]{Research partially supported by NSA grant MDA 904-00-1-0046 and NSF grant DMS-0139287. This manuscript is submitted for publication with the understanding that the United States government is authorized to reproduce and distribute reprints.}  
\markboth{$A_6$ extensions and mod $p$ Cohomology}{A.\ Ash, D.\ Pollack, and W.\ Sinnott}  
\section*{Introduction}
    In this paper we give additional computational evidence for the  generalization of Serre's Conjecture [Serre 87] proposed in [Ash-Sinnott  00] and extended in [Ash-Doud-Pollack 02].  We also propose a further  refinement --- corresponding precisely to the ``peu ramifi\'ee -- tr\`es  ramifi\'ee'' distinction in Serre --- which removes an ambiguity in the  prediction of the earlier conjecture.  
    
The conjecture sets out a precise relationship between mod $p$ Galois  representations and the mod $p$ cohomology of congruence subgroups of  $SL_n(\ZZ)$: see below for the statement in the case of 3-dimensional irreducible representations, and see [Ash-Doud-Pollack 02] for the general situation of $n$-dimensional, possibly reducible, representations. 

Let $p$ be a prime number, and let $\bFp$ be an algebraic closure of the finite field $\FF_p$ with $p$ elements. We fix an algebraic closure $\bar\QQ$ of the rational numbers $\QQ$, and let $G_\QQ=\Gal(\bar \QQ / \QQ)$ be its Galois group. We consider in this paper 3-dimensional representations of $G_\QQ$ over $\bFp$ of ``type $A_6$'': continuous irreducible Galois representations $\rho:G_\QQ\to GL_3(\bFp)$ whose image in $PGL_3(\bFp)$ is isomorphic to $A_6$.  The image of $\rho$ has larger order in these examples than any on which our conjecture been tested so far. 

We looked at twelve $A_6$ extensions of $\QQ$, taken from the tables of John W.\ Jones ( \emph{Tables of Number Fields with Prescribed Ramification: Sextics}[Jones 98]). These are in fact all $A_6$ extensions ramified at at most two primes $\le 19$. We found six of these on which it was feasible to check the conjecture. In all of these cases, the conjecture appears to be true, in the sense that a Hecke eigenclass exists in the predicted cohomology group satisfying the equality of its Hecke polynomial at $\ell$ with the characteristic polynomial of Frobenius at $\ell$ for all unramified $\ell \le 47$. 

Part of the task of checking our conjecture requires studying the local behavior of these fields: determining local Galois groups and higher ramification. These computations were initially being done entirely with PARI/GP [GP], but while we were investigating these examples, the \emph{Database of Local Fields} created by John W.\ Jones and David P.\ Roberts [Jones-Roberts 03] (with its enormously useful local fields calculator) came online. It allowed us to double-check the local field calculations we had made, and it made the remaining calculations much easier to complete. The other part of checking our conjecture requires determining the action of the Hecke algebra on the cohomology of certain congruence subgroups of $ SL_3(\ZZ)$ with mod $p$ coefficients. For this the programs developed for [Ash-Doud-Pollack 02] (which have now been extended to allow $p=2$) were used. 

We would like to thank  Koichiro Harada, Ron Solomon, John Swallow, and David Roberts  for helpful conversations and information.

\section{The Conjecture}

In this section we review and slightly refine the conjecture on mod $p$ Galois representations found in section 3 of [Ash-Doud-Pollack 02]. However, we only state the conjecture for irreducible, three-dimensional representations.  The refinement consists of resolving the ambiguity inherent in the ``prime" notation for the weight, using the concept of peu vs.\ tr\`es ramifi\'ee.  We simply copy Serre in his original conjecture for $GL(2)$ in this regard.

We choose for each prime $q$ a Frobenius element $\Frob_q$ in $G_\QQ$ and we fix a decomposition group $G_q$ of $q$ with its filtration by its ramification subgroups $G_{q,i}$ for $i\ge 0$, numbered so that $G_{q,0}$ is the full inertia group of $q$ in $G_q$.  We sometimes denote $G_{q,0}$ by $I_q$. The quotient $I_q/G_{q,1}$ is called the
\emph{tame inertia group} of $q$. We also fix a complex conjugation $\Frob_\infty \in G_\QQ$.

Start with a continuous irreducible representation $\rho: G_\QQ \to GL_3(\bFp)$.  We will define three invariants associated to $\rho$: a level, a nebentype character, and a collection of weights.

\subsection*{Level:} 
For any prime $q \ne p$ set 
\begin{equation} 
n_q =\sum_{i=0}^\infty (|G_i|/|G_0|)\  \text {dim} M/M^{G_i} ,
\end{equation}
where $G_i$ denotes the image under $\rho$ of the $i$-th higher ramification subgroup at $q$ and $M$ denotes the vector space ${\bFp}^{\,3}$ on which $\rho$ acts.  We define the level of $\rho$ to be $N =\prod_{q \ne p} q^{n_q}$. 

\subsection*{Nebentype:} 
Factor det$\rho = \eps \omega^k$, where $\omega$ is the cyclotomic character modulo $p$ and $\eps$ is a character of $G_\QQ$ unramified at $p$.  The conductor of $\eps$ divides $N$, so we can consider it as a Dirichlet character $(\ZZ/N\ZZ)^\times \to \bFp^\times$.  This is the nebentype character of $\rho$.

\subsection*{Weight:} 
This is a triple of integers $(a,b,c)$ satisfying $0 \le a-b, b-c \le p-1$ and $0 \le c \le p-2$.  We call such a triple ``$p$-restricted".  Given such a triple, there is associated an irreducible $\FF_p[GL_3(\ZZ/p\ZZ)]$-module, which is denoted by $F(a,b,c)$.  See section 2.3 of [Ash-Doud-Pollack 02].  Associated to $\rho$ there is a collection of weights, determined as follows.

Consider $\rho$ restricted to $I_p$, the inertia subgroup of $G_\QQ$ at $p$.  This can be upper-triangularized over $\bFp$, and we obtain three characters down the diagonal.  There are three cases, which we refer to by the term ``niveau".  Recall that a character of $I_p$ into $\bFp^\times$ factors through the tame inertia group $I_p/G_{p,1}$ and is called ``niveau $k$" if its image lies in $\FF_{p^k}^\times$  but in no $\FF_{p^r}^\times$  for $r < k$.

\emph{Niveau 1:} in this case all three characters are niveau 1, i.e. powers of $\omega$.  Let those powers be $A,B,C$ read from the upper left corner to the bottom right corner.  Let $(a,b,c)$ be the $p$-restricted triple congruent to $(A-2,B-1,C)$ modulo $p-1$. (We sometimes use the ``prime'' notation, and write $(a,b,c)=(A-2,B-1,C)'$.) We say that $(a,b,c)$ is one of the weights associated to $\rho$. Doing this for all different ways of upper-triangularizing $\rho|{I_p}$, we obtain the set of weights associated to $\rho$. 

\emph{Niveau 2:} Let $\psi$ be the fundamental character of niveau 2.  In this case, one of the characters is niveau 2, say $\psi^m$, another of the characters is $\psi^{pm}$, and the third character is niveau 1, say $\omega^k$.  Write $m = r+sp$ with $0 \le r-s \le p-1$ (if possible).  Let $(A,B,C)$ be either $(k,r,s)$, $(r,k,s)$ or $(r,s,k)$. Let $(a,b,c)$ be the $p$-restricted triple congruent to $(A-2,B-1,C)$ modulo $p-1$.  We say that $(a,b,c)$ is one of the weights associated to $\rho$. Doing this for all possible choices, we obtain the set of weights associated to $\rho$.

\emph{Niveau 3:} Let $\theta$ be the fundamental character of niveau 3.  In this case, one of the characters is niveau 3, say $\theta^m$, another of the characters is $\theta^{pm}$, and the third character is $\theta^{p^2m}$.  Write $m = r+sp+tp^2$ with $0 \le r-t, s-t \le p-1$ (if possible).  Let $(A,B,C)$ be the rearrangement of $(r,s,t)$ into descending order. (Note that conditions  on $r$, $s$, and $t$ imply that $C=t$, but $r$ and $s$ may not be in descending order.) Let $(a,b,c)$ be the $p$-restricted triple congruent to $(A-2,B-1,C)$ modulo $p-1$.  We say that $(a,b,c)$ is one of the weights associated to $\rho$. Doing this for all possible choices, we obtain the set of weights associated to $\rho$.

If $A-B-1$ or $B-C-1$ or both are divisible by $p-1$ there is an ambiguity in our definition of $(a,b,c)$. In niveau 1 we resolve it as follows: suppose first that $A-B-1$ is divisible by $p-1$.  We view the upper left 2-by-2 block as a representation of $\rho|{I_p}$ into $GL(2,\bFp)$.  This representation is either peu or tr\`es ramifi\'ee, in the sense of section 2.4 (ii) of [Serre 87].  If it is peu ramifi\'ee, we allow both $a=b$ and $a=b+p-1$. If it is tr\`es ramifi\'ee,  we set $a=b+p-1$.  We resolve the ambiguity if $B-C-1$ is divisible by $p-1$ in a similar way.

In the niveaux 2 and 3 we expect that any ambiguity should be resolved by allowing both possibilities, i.e., if  $A-B-1$ is divisible by $p-1$ we should allow both $a=b$ and $a=b+p-1$, and if $B-C-1$ is divisible by $p-1$ we should allow both $b=c$ and $b=c+p-1$. See the upcoming paper of Ash, Pollack and Soares [Ash-Pollack-Soares] for examples of this behavior when $p=2$.

Before we state the conjecture we briefly review the cohomology of certain subgroups of $GL_3(\ZZ)$ as Hecke modules.

For any level $N$ prime to $p$, set $\Gamma_0(N)$ to be the subgroup of $SL_3(\ZZ)$ consisting of matrices whose first row is congruent to $(*,0,0)$ modulo $N$.  If $\eps: \ZZ/N\ZZ \to \bFp^\times$ is a character, we pull it back to a character of $\Gamma_0(N)$.  If $V$ is any $\bFp[SL_3(\ZZ/p\ZZ)]$-module, pull it back to a $\Gamma_0(N)$-module $V^*$.  We denote by $V_\eps$ the $\Gamma_0(N)$-module $V^* \otimes \eps$.

The Hecke algebra $H$ acts on $H^*(\Gamma_0(N), V_\eps)$.  Let $x$ be a Hecke eigenclass and $\ell$ a prime not dividing $pN$.  Write $T_{\ell,k}x = a_{\ell,k}x$ where $T_{\ell,k}$ is the Hecke operator corresponding to the double coset $\Gamma_0(N) \text{diag}(1,\dots,\ell) \Gamma_0(N)$ with $k$ $\ell$'s, and $a_{\ell,k} \in \bFp^\times$.

If $\rho: G_\QQ \to GL_3(\bar \FF_p)$ is a continuous irreducible representation unramified outside $pN$, we say that $\rho$ is attached to $x$ if
\begin{equation}
\sum_{k=0}^3 (-1)^k \ell^{k(k-1)/2} a_{\ell,k} t^k = \text{det}(I - \rho(\Frob_\ell)t)
\end{equation} 
for all $\ell$ not dividing $pN$.

We are ready to state the conjecture.

\smallskip

\textbf{Conjecture.}\textit{ Let $\rho: G_\QQ \to GL_3(\bar \FF_p)$ be a continuous irreducible representation.  If $p > 2$, we assume that $\rho(\Frob_\infty)$ has eigenvalues $1,1,-1$ or $1,-1,-1$.  Let $N$ be the level and $\eps$ the nebentype character associated to $\rho$. Then for any weight $(a,b,c)$ associated to $\rho$, there exists a Hecke eigenclass $x$ in $H^3(\Gamma_0(N), F(a,b,c)_\eps)$ with $\rho$ attached.}

\section{Fields with Galois group $A_6$ ramified at two small primes}

We used the tables of Jones [Jones 98] to obtain a complete list of all sextic extensions of $\QQ$ which are at most ramified at two primes $\le 19$ and whose Galois closure has Galois group $A_6$.  There are 24 such fields in these tables, each ramified at two primes $\le 19$ (and at infinity).  However, $A_6$ has two conjugacy classes of subgroups of index 6: one is represented by the ``natural'' subgroups isomorphic to $A_5$, i.e.\ the stabilizers of the elements $1,\ldots, 6$ under the natural permutation representation of $A_6$, and the two classes are interchanged by an outer automorphism of $S_6$.  Hence sextic fields whose Galois closure has group $A_6$ will occur in pairs with the same Galois closure.   

We can determine the pairing as follows.  Since an outer automorphism of $S_6$ switches the two conjugacy classes of elements of order 3 (the 3-cycles and the double 3-cycles) and preserves the other $S_6$-conjugacy classes in $A_6$, we can identify sextic fields with the same Galois closure by examining the cycle structure of the Frobenius automorphism $\Frob_q$ for various primes $q$. In more detail: let $K/\QQ$ be a Galois extension of \QQ\ with Galois group $G\simeq A_6$, let $H$ and $H'$ be representatives of the two conjugacy classes of subgroups of $G$ of index 6. Then $K^H$ and $K^{H'}$ will be non-conjugate sextic extensions of \QQ\ whose Galois closure is $K$. If we let $G$ act on the cosets of $H$ and $H'$, we obtain two homomorphisms $\phi, \phi':G\to S_6$ (depending on a numbering of the cosets): both homomorphisms give an isomorphism of $G$ with $A_6$, and we have $\phi'=\alpha\circ\phi$, for some outer automorphism $\alpha$ of $S_6$. If $\Frob_q\in G$ is a Frobenius for $q$, then the cycle structures of $\phi(\Frob_q)$ and $\phi'(\Frob_q)$ are the same unless $\Frob_q$ has order 3 ($\alpha$ switches the 3-cycles and the double 3-cycles, but preserves the other $S_6$-conjugacy classes in $A_6$). 

Thus there are twelve Galois extensions with Galois group $A_6$ arising from Jones's tables. These extensions are listed in Table 1. Since Jones's tables are complete, this is a \emph{complete} list of all $A_6$ extensions of $\QQ$ at most ramified at two primes $\le 19$.

\renewcommand{\arraystretch}{1.1}
\begin{table}[hbt]
\begin{tabular}{l|l}
location & polynomial  \\
    \hline
 2,3: \#55 or \#56 &$x^6+3x^5+3x^4 +2x^3 -3x^2 -3x -1$ \\
 2,3: \#57 or \#60 &$x^6+8x^3+9x^2-6$ \\
 2,3: \#58 or \#61 &$x^6+6x^4-4x^3-3x^2-12x-12$  \\
 2,3: \#59 or \#62 &$x^6-12x^3+21x^2+12x-34$  \\
 2,5: \#17 or \#18 &$x^6-2x^5+15x^4+50x^2-4x-82$\\
 3,5: \#7 or \#10  &$x^6-5x^3+45x^2-99x-15$  \\
 3,5: \#8 or \#9   &$x^6+3x^5+15x^4+25x^3+45x+60$  \\
 3,7: \#3 or \#4   &$x^6+3x^5+3x^4-9x^3-18x^2+9x+18$  \\
 3,13: \#9 or \#10 &$x^6+3x^5+3x^4+2x^3+3$ \\
 3,19: \#3 or \#4  &$x^6-3x^5-3x^4+14x^3-12x+9$\\
 5,17: \#1 or\#2    &$x^6-2x^5+5x^2-11x-13$  \\
 13,19: \#1 or\#2   &$x^6-4x^4-15x^3-15x^2-8x+4$
\end{tabular}
\medskip
\caption{$A_6$ Sextics}
\end{table}
Under ``location'' we have indicated where in Jones's tables each $A_6$ extension can be found; e.g.\ ``2,3:\#55 or \#56'' means that entries \#55 and \#56 in the table \emph{New sextic fields ramifying above \{2,3\}} are non-conjugate sextic fields with the same $A_6$ extension as their normal closure. The polynomial listed in each case generates the first field of the pair. Henceforth we refer to these $A_6$ extensions by the first sextic in each pair.

We used this table to find representations of type $A_6$ on which we could test our conjecture. Let $K$ be the splitting field of one of these polynomials. Then $\Gal(K/\QQ)\simeq A_6$. If $p\ne 3$, $A_6$ has no irreducible 3-dimensional representations over $\bFp$, but its triple cover $3.A_6$ does, so when $p\ne 3$ we have to determine whether $K$ is contained in a $3.A_6$ extension $\tK/\QQ$. If such an extension $\tK$ exists, we obtain from it a representation $\rho:G_\QQ\to GL_3(\bFp)$ (via $G_\QQ\to \Gal(\tK/\QQ)\to GL_3(\bFp))$. On the other hand, if we take $p=3$, $A_6$ does have representations into $GL_3(\bF_3)$, so we can use $K$ directly to obtain a representation $\rho:G_\QQ\to GL_3(\bF_3)$. (We discuss the 
group $3.A_6$  and its representations more fully in the next section.)  Since our conjecture is stable under twisting by Dirichlet characters,  we can suppose without loss of generality that (when $p\ne 3$) $\tK$ is ramified at the same primes as $K$. This can be seen as follows, using a technique due to Tate ([Serre 77],\S6). Suppose that $\ell$ is a rational prime which is ramified in $\tK$ but not in $K$. Then $\rho(I_\ell)$ equals the center $Z$ of $SL_3(\bFp)$. Let $G_\ell$ be a decomposition group of $\ell$ in $G_\QQ$: since $\rho(G_\ell)/Z$ is isomorphic to the decomposition group of $\ell$ in $\Gal(K/\QQ)$, it is cyclic, and therefore $\rho(G_\ell)$ is abelian. Therefore the commutator subgroup $G_\ell'$ of $G_\ell$ is contained in $\ker\rho$, and $G_\ell'\subseteq I_\ell$; hence $\rho|I_\ell$ factors through $I_\ell/G_\ell'$. Now  $I_\ell/G_\ell'$ may be identified with $\Gal(\QQ_\ell(\mu_{\ell^{\infty}})/\QQ_\ell)$,  which in turn can be identified with $\Gal(\QQ(\mu_{\ell^{\infty}})/\QQ)$. (Here $\mu_{\ell^{\infty}}$ is the group of roots of unity of order a power of $\ell$.) Hence there is a cubic character $\chi_\ell$ of $G_{\QQ}$,  unramified outside $\ell$, such that
\[
\rho(\tau)=\chi_\ell(\tau)I, \text{ for $\tau\in I_\ell$.}
\]
Here $I$ is the $3\times 3$ identity matrix. So $(\prod_\ell\chi_\ell)^{-1}\cdot\rho$ is only ramified at the primes that ramify in $K$. (Here $\ell$ runs over rational primes that ramify in $\tK$ but not in $K$.) Thus we may assume that the level $N$ of $\rho$ is divisible only by the primes (other than $p$) which ramify in $K$. Note that $\rho$ still has determinant 1, so it has trivial nebentype, and still has image $3.A_6$ in $SL_3(\bFp)$. Sometimes the level $N$ can be lowered further by twisting by a suitable Dirichlet character ramified at a prime ($\ne p$) that does ramify in $K$; this can change the nebentype.  

There are two aspects to testing our conjecture:
\begin{itemize}
     \item We have to study the field $\tK$ (or just $K$, when $p=3$), to determine the level $N$ and weights $(a,b,c)$ attached to $\rho$; and we need to determine the characteristic polynomials of $\rho(\Frob_\ell)$ for $\ell\nmid pN$.  The main difficulty here is the analysis of the ramification groups, which is required to determine the level and the weight. We give examples of these calculations below.
     \item Once the level $N$, nebentype $\eps$, and weights $(a,b,c)$ have been found, we have to look for eigenclasses $x$ for the action of the Hecke algebra on $H^3(\Gamma_0(N), F(a,b,c)_\eps)$ which appear correspond to $\rho$. This means checking that the equation $(2)$ holds for primes $\ell\nmid pN$, $\ell \le 47$.  For these calculations, we use the programs developed for [Ash-Doud-Pollack 02] (which have now been extended to allow $p=2$); for a description of these programs, see [Ash-Doud-Pollack 02]. The main difficulties that arise here are limitations on the feasibility of the computations if the level $N$ or the weight are too large. For this reason, we have always taken $p$ to be one of the primes that ramify in $K$, since this substantially reduces the level.
     \end{itemize}

\section{modular representations of $3.A_6$}

The Schur multiplier of $A_6$ has order 6, so the
universal central extension of $A_6$ is a six-fold cover of
$A_6$, denoted by $6.A_6$, and the quotient of $6.A_6$ by
its central subgroup of order 2 is a triple cover of $A_6$
denoted by $3.A_6$. $3.A_6$ has faithful 3-dimensional
complex representations; explicit generators for the
image of one of these representations can be found in 
[Crespo-Hajto]. 

If $p\ne 3$, $A_6$ has no 3-dimensional linear representations over $\bFp$, but it does have representations into $PGL_3(\bFp)$. 
The map $SL_3(\bFp)\to PGL_3(\bFp)$ is surjective, and its kernel
has order 3, and therefore  any representation of $A_6$ into $PGL_3(\bFp)$ can be lifted uniquely to a representation of $3.A_6$ into $SL_3(\bFp)$. On the other hand, $A_6$ does have 3-dimensional linear representations over $\bF_3$. 

All of these representations can be obtained by reduction mod $p$ from the 3-dimensional complex representations of $3.A_6$.  There are 4 such representations, which can be realized over $F=\QQ(\zeta_3, \sqrt{5})$.  In fact, the whole character table of $3.A_6$ lies in $F$, and all the representations of $3.A_6$ can be realized over $F$: for any finite group $G$, the Schur indices of $F[G]$ divide the order of the center of $G$, and also divide $\phi(m)$, where $m$ is the exponent of $G$ (see for example Serre [Serre 71], pg.\ 109).  Since the center of $3.A_6$ has order 3 and the exponent of $3.A_6$ is 60 (and $\phi(60)=16$), the Schur indices of $F[3.A_6]$ are all 1.  So all the complex representations of $3.A_6$ can be realized over $F$. The four 3-dimensional representations of $3.A_6$ are conjugate over $\QQ$.

Let $p$ be a rational prime, let \p\ be a prime of $F$ lying above $p$, let $\tilde\rho: 3.A_6\to GL_3(\o_\p)$ be a realization of one of the 3-dimensional representations of $3.A_6$ over the integers $\o_\p$ of $F_\p$, and let $\rho$ be the reduction of $\tilde\rho$ mod \p.  All the 3-dimensional modular representations of $3.A_6$ (up to similarity) arise in this way; this can be seen from the decomposition matrices of $3.A_6$, which are available from [GAP], for example.

We can summarize 3-dimensional modular representations of $3.A_6$ as follows:
\begin{itemize}
    \item $p>5, p\equiv 1,4 \bmod 15$: There are 4 irreducible 
    3-dimensional representations of $3.A_6$, each defined over 
    $\FF_p$.
    \item$p>5, p\not\equiv 1,4 \bmod 15$: There are 4 irreducible 
    3-dimensional representations of $3.A_6$, each defined over 
    $\FF_{p^2}$. They are conjugate in pairs over $\FF_p$.
    \item $p=5$: There are 2 irreducible 
    3-dimensional representations of $3.A_6$, each defined over 
    $\FF_{25}$. They are conjugate over $\FF_{5}$.
    \item $p=3$: There are 2 irreducible 
    3-dimensional representations of $3.A_6$, each defined over 
    $\FF_{9}$. They are conjugate over $\FF_{3}$.
    \item $p=2$: There are 4 irreducible 
    3-dimensional representations of $3.A_6$, each defined over 
    $\FF_{4}$. They are conjugate in pairs over $\FF_2$.
\end{itemize}

All these representations are faithful except the mod 3 representations, for which the center lies in the kernel: these are 3-dimensional representations of $A_6$ in $GL_3(\FF_9)$. 

\section{mod $p$ Galois representations of type $A_6$}

In this section we describe the mod $p$ Galois representations of type $A_6$ that arise from the twelve $A_6$-extensions listed in Table 1, and give examples of the computation of the level $N$ and the associated weights, and the characteristic polynomials of Frobenius. As explained in \S 2, we always assume that $p$ is one of the ramified primes. We discuss the cases with $p=3$ separately since we don't need to worry about the existence of a triple cover in those cases.

\subsection*{mod 3 Galois representations of type $A_6$ }

Let $K$ be one of the nine $A_6$ fields in Table 1 that are ramified at 3. We take $p=3$ and let $\rho$ be a  representation  of $G_\QQ$ into $GL_3(\bF_3)$ which cuts out $K$, i.e. $K=\bar\QQ^{\ker\rho}$. As discussed above, there are two such representations: we may assume that the image of $\rho$ lies in $GL_3(\FF_9)$, and the other possible choice for $\rho$ is $\rho^\phi$, where $\phi$ is the nontrivial element of $\Gal(\FF_9/\FF_3)$. $\rho$ has trivial nebentype.

For three of the nine fields $K$, we can lower the level by twisting by a suitable quadratic Dirichlet character $\eps$; for the remaining fields, we let $\eps=1$. Then for all nine cases the representation is $\rho\eps$ and the nebentype is $\eps$. In Table 2 below, we list the nine fields which are ramified at 3, the level $N$ of the representation $\rho\eps$, and the nebentype $\eps$. (We omit the nebentype $\eps$ when it is trivial.) $\omega_4$ denotes the Dirichlet character corresponding to $\QQ(i)$ and $\psi_8$ the Dirichlet character corresponding to $\QQ(\sqrt 2)$, considered as taking values in $\FF_3$.

\begin{table}[hbt]
\begin{tabular}{l|l}
$K$ & $N$, $\eps$\\
  \hline
 2,3: \#55 & $2^8$\\
 2,3: \#57 & $2^{7}$, $\psi_8$\\
 2,3: \#58 & $2^{7}$, $\omega_4\psi_8$\\
 2,3: \#59 & $2^{8}$, $\omega_4$ \\
 3,5: \#7  & $5^4$\\
 3,5: \#8  & $5^4$\\
 3,7: \#3  & $7^2$\\
 3,13: \#9 & $13^2$\\  
 3,19: \#3 & $19^2$\\
\end{tabular}
\medskip
\caption{Levels of mod 3 representations}
\end{table}        
We illustrate the computation of the level for the first field above (2,3:\#55).  Let $K$ be the splitting field of $T(x)=x^6+3x^5+3x^4 +2x^3 -3x^2 -3x -1$, $G=\Gal(K/\QQ)\simeq A_6$. We need to find the images of the ramification groups $G_{2,i}$ of 2 in $G$.  To this end we let $L$ be the splitting field of $T$ over $\QQ_2$, and let $D=\Gal(L/\QQ_2)$; $D$ may be identified with a decomposition group of 2 in $G$. Let $G_i=\rho(G_{2,i})$; the $G_i$'s may be identified with the ramification groups of $L/\QQ_2$.

The following computations are easily done in PARI/GP. If we factor $T$ over $\QQ_2$ we find that $T=H\cdot Q$, where $H$ is an irreducible quartic and $Q$ is an irreducible quadratic. Using the results of the Appendix, we find that $g(x)=x^4 + 4x^3 -14x^2 + 3$ generates the same extension of $\QQ_2$ as $H$, and $q(x)=x^2+x+1$ generates the same extension of $\QQ_2$ as $Q$.  Note that $g$ and $q$ are both irreducible over $\QQ_2$: $g(x-1)=x^4 - 20x^2 + 36x - 14$ is an Eisenstein polynomial, and $q$ is irreducible mod 2.

Let $g_1(x)=g(x-1)=x^4 - 20x^2 + 36x - 14$, and let $h_1$ be the resolvent cubic of $g_1$: $h_1=x^3 + 20x^2 + 56x - 176$.  Then $h_2(x)=h_1(2x)/8= x^3 + 10x^2 + 14x - 22$ is Eisenstein, so $h_1$ is irreducible; the discriminant of $h_1$ has the form $2^8\cdot u$, where $u\equiv 5\bmod 8$---which is not a square, so the Galois group (over $\QQ_2$) of $h_1$ is $S_3$, and the Galois group (over $\QQ_2$) of $g_1$ is $S_4$. 

The splitting field of $g_1(x)$ over $\QQ_2$ contains the discriminant field of $g_1$, namely $\QQ_2(\sqrt5)$, which is the unramified quadratic extension of $\QQ_2$, and is also the splitting field of $q$.  Hence the splitting field of $g_1$ is also the splitting field $L$ of $T(x)$ over $\QQ_2$. Thus $D\simeq S_4$. In what follows we will identify $D$ with $S_4$.

The ramification index of $L/\QQ_2$ is divisible by 4 (from $g_1$) and 3 (from $h_2$), so equals 12 or 24.  But it can't be 24, since $\QQ_2(\sqrt5)/\QQ_2$ is unramified.  So $G_0=A_4$ and $G_1=V_4$ (since $G_1$ is the Sylow 2-subgroup of $G_0$). Since $S_4$ has no normal subgroups between $V_4$ and 1, we have $G_1=G_2=\ldots=G_r$ and $G_{r+1}=1$ for some $r\ge 1$; we need to determine $r$.

To study the higher ramification, let $\alpha$ be a root of $g_1(x)$ in $L$, $\beta$ a root of $h_2(x)$, and let $\pi=\beta/\alpha$: then $\ord_2(\pi)=\frac{1}{3}-\frac{1}{4}=\frac{1}{12}$, so $\pi$ is a local parameter in $L$ (here $\ord_2$ is normalized so that $\ord_2(2)=1$). Let $\ord_L$ be the valuation on $L$, normalized so that $\ord_L(\pi)=1$. Let $\sigma\in V_4$. Then $\sigma\beta=\beta$ (because $\beta$ has degree 3 over $\QQ_2$) so that we have
\begin{align*}
\ord_L(\sigma\pi-\pi)&=
\ord_L\left(\beta\frac{\alpha-\sigma\alpha}{\alpha\,\sigma\alpha}\right)\\
&=4-3-3+\ord_L(\alpha-\sigma\alpha)\\
	 &=\ord_L(\alpha-\sigma\alpha)-2.
\end{align*}
Let $\nu=\ord_L(\sigma\pi-\pi)$ for $\sigma\ne1, \sigma\in V_4$; $\nu$ doesn't depend on $\sigma$ because there are no ramification groups between $V_4$ and 1. Since $g_1'(\alpha)=\prod_{\sigma\in V_4, \sigma\not=1}(\alpha-\sigma\alpha)$ and $\norm_{L/\QQ_2}(g_1'(\alpha))=disc(g_1)^6$, we get
\begin{align*}    
\ord_L \norm_{L/\QQ_2}(g_1'(\alpha))&=\ord_L(disc(g_1)^6)\\
 24\cdot3\cdot(\nu+2)&=6\cdot8\cdot12;
\end{align*}
so $\nu=6$, and so $r=5$.

It follows that $G_i=\rho(G_{2,i})$ contains a subgroup of $GL_3(\FF_9)$ isomorphic to $V_4$ for $i=0\ldots 5$.  We can diagonalize this subgroup and assume that $G_i$ (for $i=0\ldots5$) contains the subgroup of diagonal matrices of the form $diag(\pm1, \pm1, \pm1)$ (with determinant 1).  It follows that $M^{G_i}=0$ for such $i$.  So the power of 2 in the level of $\rho$ is 8, as can be computed from the definition (1) in $\S2$. So $N(\rho)=2^8$. 

To compute the weight(s) attached to $\rho$, we find the possible upper-triangularizations of $\rho(I_p)$ (with $p=3$). Let $S$ be the splitting field of $T$ over $\QQ_3$: we need to determine the inertia group of $S/\QQ_3$ (which may be identified with $\rho(I_3)$). For this we can use the ``local field calculator'' attached to the local field database of Jones and Roberts [Jones-Roberts 03], which identifies the field $S$ as an extension of $\QQ_3$ with inertia group $C_3\times C_3$. Hence the elements of $\rho(I_3)$ are unipotent in $GL_3(\bF_3)$. It follows that there is just one weight attached to $\rho$ by our conjecture, namely $(-2,-1,0)'=(2,1,0)$.

The difficulties in computing $N(\rho)$ in this example arise from the fact that 2 is wildly ramified in $K$. For the examples ramified at 3 and $q$ with $q\ne 2, 5$, the ramification  will be tame, and $N(\rho)$ is  easier to determine. For example, consider the field ramified at 3 and 13  (3-13 \#9). The polynomial is $T(x)=x^6+3x^5+3x^4+2x^3+3$, which factors over $\QQ_{13}$ into two cubics $g_1$ and $g_2$, both of which have square discriminants (in $\QQ_{13}$) and therefore cyclic Galois groups. $g_1(x-5)$ is Eisenstein, while $g_2$ is irreducible mod 13. So the splitting field of $T$ over $\QQ_{13}$ has Galois group $D\simeq C_3\times C_3$; the inertia group $G_0$ has order 3, and we have $G_1=1$, since the ramification is tame. 

We now need to know the dimension $M^{G_0}$: since $G_0$ is generated by a 3-cycle $\tau$ in $A_6$ (in terms of its action on the roots of $T$), and since a 3-cycle in $A_6$ normalizes a $V_4$ (i.e., we may view the 3-cycle as an element of $A_4\subseteq A_6$) it follows that we may suppose that $\rho(\tau)$ is a permutation  matrix of order 3 in $GL_3(\FF_9)$, which fixes a 1-dimensional subspace of $M$. From this it is easy to see from the definition (1) that the level of $\rho$ is $13^2$.

\subsection*{mod $p$ Galois representations of type $A_6$, $p\ne3$}

Let $K$ be a Galois extension of $\QQ$ with Galois group isomorphic to $A_6$, and let $p\ne3$.  There is always a projective representation $\tilde\rho: G_\QQ\to PGL_3(\bFp)$ which cuts out $K$, and we need to examine the question of whether $\tilde\rho$ lifts to a special linear representation $\rho: G_\QQ\to SL_3(\bFp)$.  (Note that $PSL_3(\bFp)=PGL_3(\bFp)$.)  Now, a theorem of Neukirch (see  [Neukirch 73], Satz 2.2) says that such a $\rho$ exists precisely if the ``local lifting problem'' can be solved for each place $v$ on $\QQ$.  This means that for each place $v$ of $\QQ$, we have to find a lifting
\[
\rho_v: G_v\to SL_3(\bFp)
\]
of $\tilde\rho\vert{G_v}:G_v \to PGL_3(\bFp)$. Here $G_v$ denotes the decomposition group at $v$. Also, the local lifting problem at $v$ is always solvable if the order of $\tilde\rho(G_v)$ is not divisible by 9. This is an observation of Feit's (see [Feit 89], \S6).  Indeed, let $Y=\tilde\rho(G_v)$, and suppose that $Y$ does not have order divisible by 9. Let $H$ be the inverse image of $Y$ in $SL_3(\bFp)$. Then $H$ contains the center $Z\simeq C_3$ of $SL_3(\bFp)$. Let $W$ be the  Sylow 3-subgroup of $Y$, and $X$ the Sylow 3-subgroup of $H$. Then either $W=1$ and  $X=Z$ or $W\simeq C_3$ and $X\simeq Z\times C_3$ (since $3.A_6$ has no elements of order 9). So the class of $X$ in $H^2(W, Z)$ is trivial, and therefore the class of $H$ in $H^2(Y,Z)$ is trivial ($\text{res}_W^Y:  H^2(Y,Z)\to H^2(W,Z)$ is injective). Therefore $H$ splits as $Z\times Y$. So the local lifting problem is solvable. The theorem of Neukirch and the observation of Feit were pointed out to us by John Swallow.

Suppose that $v=q$ is unramified in $K$: if $g\in SL_3(\bFp)$ is any lifting of $\tilde\rho(\Frob_q)$, we can define a lifting $ \rho_q$ by setting $\rho_q(\Frob_q)=g$ and $\rho_q(I_q)=1$ (the point being that $G_q/I_q\simeq \hat \ZZ$). If $v=\infty$, then 9 does not divide the order of $\tilde\rho(G_\infty)$. Hence the possibly delicate cases occur when $v$ is a finite prime which ramifies in $K$. 

Table 3 below sets out  what we have been able to find out about the question of whether the fields in Table 1 embed in a $3.A_6$ extension. We would like to thank David Roberts for resolving for us the harder cases (2,3\#55, 3,5\#8, and 3,7\#7).  Roberts also finds eighteenth degree polynomials whose splitting fields give the $3.A_6$ extensions when they exist. For a discussion of these questions, see his upcoming paper [Roberts].  

\renewcommand{\arraystretch}{1.1}
\begin{table}[hbt]
\begin{tabular}{l|l}
$A_6$ field & $3.A_6$?  \\
    \hline 
 2,3: \#55 &yes  \\	 
 2,3: \#57 &yes  \\	 
 2,3: \#58 &yes  \\	 
 2,3: \#59 &yes  \\	
 2,5: \#17 &yes\\
 3,5: \#7  &yes\\
 3,5: \#8  &yes\\
 3,7: \#3  &yes\\
 3,13: \#9 &no \\
 3,19: \#3 &yes\\
 5,17: \#1 &yes\\  
 13,19: \#1 &no
\end{tabular}
\medskip
\caption{$3.A_6$ extensions}
\end{table}  

\subsection*{An example: the $A_6$ extension ramified at 5 and 17}

We discuss this example in some detail.  (The calculations were done with PARI/GP, v.  2.1.4.)  Let $K$ be the splitting field of $T(x)=x^6-2x^5+5x^2-11x-13$.  This polynomial is 5,17:\#1 from Table 1. The discriminant of the sextic field generated by a root of $T$ is $5^8 17^2$, and $\Gal(K/\QQ)\simeq A_6$.

First we show that $K$ is contained in a $3.A_6$ extension $\tilde K$ of $\QQ$. Let $\tilde\rho:G_\QQ\to PGL_3(\bFp)$ be the projective representation that cuts out $K$, where $p$ is any prime except 3 (later we will take $p=5$).
\begin{itemize}
    \item Over $\QQ_5$, $T$ factors into a linear factor and a totally ramified quintic.  So the Galois group of $T$ over $\QQ_5$ is a subgroup of $A_5$, hence does not have order divisible by 9.      
    \item Over $\QQ_{17}$, $T$ factors into three quadratics, one unramified, one equal to $\QQ_{17}(\sqrt{17})$; the third lies in the compositum of the other two. So the Galois group of $T$ over $\QQ_{17}$ is a four-group, and 4 is not divisible by 9.
\end{itemize}

Hence the local lifting problems are all solvable, so there is a representation $\rho:G_\QQ\to SL_3(\bFp)$.  The kernel of the map $SL_3(\bFp)\to PGL_3(\bFp)$ is the center of $SL_3(\bFp)$, and has order 3: so the fixed field of the kernel of $\rho$ is a $3.A_6$-extension $\tK$ of $\QQ$ containing $K$.  Let $\tG=\Gal(\tK/\QQ)\simeq 3.A_6$; then $Z=\Gal(\tK/K)$ is the center of $\tG$.

Next we show that we can twist $\rho$ by Dirichlet characters so that $\tK/K$ becomes an unramified extension.  As noted in \S2 above, we may assume that $\tK/K$ is unramified outside 5 and 17. Let $k$ be the sixth degree extension of $\QQ$ obtained from a root of $T$.  Then $\Gal(\tK/k)$ is a triple cover of $A_5$; but since the Schur multiplier of $A_5$ is 2, $\Gal(\tK/k)\simeq C_3\times A_5$.  Hence there is a cyclic cubic extension $\tk/k$ which gives $\tK$, i.e.\ $\tK = \tk K$.  Since we can and have assumed that $\tK/K$ is unramified outside 5 and 17, the conductor of $\tk/k$ must divide $5\cdot17$ (since the ramification over 5 and 17  in $\tk/k$ is must be tame).

According to PARI/GP (\texttt{bnrinit(bnfinit(T),5*17).clgp} gives us the order and the structure of the ray class group of conductor $5\cdot17$) the ray class group with conductor $5\cdot17$ has a unique quotient of order 3, which is in fact the ideal class group of $k$.  So $\tk/k$ is the Hilbert class field of $k$, and so $\tK/K$ is unramified, too.

We now take $p=5$, and find the level of $\rho:G_\QQ\to GL_3(\bF_5)$.  Let $I_{17}$ be an inertia group of 17 in $G_\QQ$; then the image $I$ of $I_{17}$ in $\tG$ has order 2. Let $s$ be an element of $I_{17}$ which maps to the nontrivial element of $I$; then we may assume that $\rho(s)= \diag(1,-1,-1)$ (since $\rho(s)$ has order 2 and $\det\rho=1$). Hence $M^{G_0}$ is 1-dimensional.  (Recall that $G_0=\rho(I_{17})\simeq I$.) On the other hand $G_1=1$, so the power of 17 in the level of $\rho$ is 2, according to (1).  Thus the level of $\rho$ is $17^2$.  

We can in fact lower this by twisting by the quadratic character $\eps_{17}$ of conductor 17. Let $\rho'=\eps_{17}\rho$. Then $\rho'(I_{17})$ still has order 2, since tame ramification is cyclic. Since $\eps_{17}(s)= -1$ (this is because $\eps_{17}$ corresponds to the local extension $\QQ_{17}(\sqrt{17})$, and $s$ restricts to the nontrivial automorphism of $\QQ_{17}(\sqrt{17})$), $\rho'(s)=\diag(-1,1,1)$, so $M^{G_0}$ is 2-dimensional for the twisted action.  So $\rho'$ has level 17 and nebentype $\eps_{17}$. 

\subsubsection*{Weights:}
To determine the weights predicted by the conjecture, we need to know how 5 ramifies in $\tK$. As noted above, $T$ factors over $\QQ_5$ into a linear factor and a totally ramified quintic $g(x)$. The polynomial $g(x+1)$ is Eisenstein, and applying Proposition 1 \emph{bis} of the Appendix, we see that the polynomial $g_1(x)=x^5-5x^4-5$ generates the same extension of $\QQ_5$. The Galois group of $g_1$ is either $C_5$ or $D_{10}$ or $M_{20}$. $M_{20}$ is ruled out because the discriminant of $g_1$ is a square. We can use the PARI/GP routine \texttt{polcompositum} to find a polynomial $h$ (of degree 20) for an extension of $\QQ$ containing two roots of $x^5-5x^4-5$. $h$ factors over $\QQ_5$ into two tenth degree polynomials; we take either of these, apply \texttt{polred}, and factor the resulting polynomials mod 5: some of these factor as $q^5 \bmod 5$, where $q$ is an unramified quadratic. This means that the Galois group of $g_1$ is $D_{10}$, and that the splitting field of $g_1$ over $\QQ_5$ is a cyclic totally ramified fifth degree extension of the unramified quadratic extension of $\QQ_5$. Thus the inertia group of 5 in $\Gal(K/\QQ)$ is cyclic of order 5, and since $\tK/K$ is unramified at 5, we can say that the inertia group of 5 in $\tG$ is also cyclic of order 5. Hence the image $\rho'(I_5)$ of inertia in $GL_3(\bF_5)$ is unipotent, so the unique weight predicted by our conjecture is $(-2,-1,0)'(\eps_{17})=(6,3,0)(\eps_{17})$. 

\subsubsection*{Characteristic polynomials of Frobenius:} 
 We used the symbolic algebra system GAP [GAP] for information about the conjugacy classes and the characters of $A_6$ and $3.A_6$, and we use the notations of that system (in particular, the labeling of conjugacy classes) in what follows.  Let $X$ be the character of the representation $\rho$; then $X$ is the reduction of the character of a 3-dimensional representation of $3.A_6$ modulo a suitable prime above 5. The character values of $X$ are as follows:
\[
\begin{array}{rrrrrrrrrrrrrr}
A6:    &1a &1a &1a &2a &2a &2a &3ab &4a &4a  &4a  &5ab &5ab  &5ab \\
3.A_6: &1a &3a &3b &2a &6a &6b &3cd &4a &12a &12b &5ab &15ac &15bd\\
X:     &3  &3z &3z'&-1 &-z &-z'&0   &1  &z   &z'  &-2  &-2z  &-2z'\\
\end{array}
\]
Here $z$ is a cube root of unity in $\FF_{25}$, $z'=z^5$ its conjugate.  Recall that there are two 3 dimensional representations of $3.A_6$ over $\bF_5$; the other one is obtained by switching $z$ and $z'$.  We have collapsed the table according to the values of $X$: thus we have written $3cd$  for the union of the two conjugacy classes $3c$ and $3d$, on which $X$ takes the  common value 0. We call these collapsed classes ($3cd$, $5ab$, etc.) ``coarse'' conjugacy classes; we only need to know in which coarse conjugacy class an element of $3.A_6$ lies to determine the value of $X$ on that element. Above each conjugacy class of $3.A_6$ we have listed its image in $A_6$. 

For each prime $\ell\ne 5, 17$, we need to determine the characteristic polynomial $\det(1-\rho'(\Frob_\ell)t)$ that appears in (2). For this it is enough to determine the coarse class of $\Frob_\ell$ in $3.A_6$, since \[
\det(1-\rho'(\Frob_\ell)t)=1-\eps_{17}(\ell)X(\Frob_\ell)t+ X(\Frob_\ell^{-1})t^2 - \eps_{17}(\ell)t^3.
\]

From the factorization of $T \bmod \ell$, we determine the cycle structure of $\Frob_\ell$ as an element of $A_6$, which determines the coarse conjugacy class of $\Frob_\ell$ in $A_6$ (the usual problem of distinguishing $5a$ and $5b$ thus disappears here).

Let $Z$ be the center of $3.A_6$, and fix a generator $c$ of $Z$. In the calculation that follows, we assume that $X(c)=3z$, i.e.\ that $c$ belongs to class $3a$. (If $X(c)=3z'$, the calculation would be similar and the result would be the same.) Suppose that $g$ is an element of $A_6$ of order prime to 3: then there is a \emph{unique} element $g'$ of $3.A_6$ of the same order which maps to $g$; and the other elements of $3.A_6$ which map to $g$ are $cg'$ and $c^2g'$.

Let $\ell$ be a prime other than 5 or 17, and let  $s=\Frob_\ell | \tK$ be a Frobenius for $\ell$ in $\tG$. We distinguish two cases.

\textsl{Case 1:} $g=s|K$ has order prime to 3; then $s=c^ig'$ for some $i$, and we would like to determine $i$. We do this as follows.

    Let $\P$ be a prime above $\ell$ in $K$, and let $\tilde\P$ be a prime above $\P$ in $\tK$; we may assume that $s$ is the  Frobenius associated to $\tilde\P/\ell$. Suppose that $f$ is the residue degree of $\P/\ell$: then
    \[
 s^f = \text{the Frobenius associated to $\tilde\P/\P$} = (\P, \tK/K),
    \] 
 so
    \[
	    c^{if} = (\P,\tK/K),
    \]
since $f$ is the order of $g$ and $g'$. So $c^{i} = (\P,\tK/K)^f$ ($f$ is its own inverse mod 3), and thus
    \[
         s = g'(\P, \tK/K)^f.
    \]     
We can ask PARI/GP to calculate $(\P, \tK/K)$ by calculating the class of $\norm_{K/k}\P$ in the ideal class group of $k$: we have 
    \[
	       (\norm_{K/k}\P, \tk/k)= (\P, \tK/K)|\tk
    \]
Finally, we need to find the coarse class of $s$. The class of $g=\Frob_\ell | K$ is determined by the cycle decomposition, and $g'$ then lies in the unique class in $3.A_6$ with the same order. $(\P, \tK/K)^f$ lies in $1a$, $3a$, or $3b$, according to whether it equals 1, $c$, or $c^2$. The class of $s=\Frob_\ell|\tK$ is then found by the rules $3a.5ab$ = $15ac$, etc. These relations can be read off from the character values of $X$.

\textsl{Case 2:} $g=s|K$ has order 3.  There is only one coarse conjugacy class whose elements have order 3 in $A_6$, namely $3ab$, and it lifts to a unique coarse conjugacy class in $3.A_6$, namely $3cd$.  So the coarse class of $s=\Frob_\ell|\tK$ is $3cd$ in this case.

%
%
%
%
%
%

We illustrate this calculation with an example. Suppose that $ \ell=2$, which has residue degree 5 in $K$; if $\P$ is a prime above 2 in $K$, then we claim that $(\P, \tK/K)=c$.

We first calculate the decomposition of 2 in $k$:

\begin{verbatim}
? T=x^6-2*x^5+5*x^2-11*x-13;
? k=nfinit(T);
? idealprimedec(k,2)
% = [[2, [1, 1, 0, 0, 0, 0]~, 1, 1, [0, 1, 0, 0, 1, 1]~], 
     [2, [2, -1, 6, -10, -7, 23]~, 1, 5, [1, 1, 0, 0, 0, 0]~]]
\end{verbatim}

So there are two primes above 2 in $k$, call them $\p_1$ and $\p_2$; $\p_2$ has residue degree 5, so that if we choose a prime $\P$ of $K$ above $\p_2$, we will have $f(\P/\p_2)=1$ and $\p_2=\norm_{K/k}(\P)$. So $\p_2$ can be used directly to find $(\P, \tK/K)$:
\begin{verbatim}
? bnfisprincipal(k,\%[2] )
% = [[1]~, [-3/5, 3/10, 2/5, -6/5, -9/10, 23/10]~, 344]
\end{verbatim}
This says that $\p_2$ lies in the ideal class of the generator that PARI/GP has chosen for the ideal class group of $k$: so $(\P, \tK/K)=c$, as claimed, since we are identifying $\Gal(\tK/K)$ with $\Gal(\tk/k)$. 

We find that the coarse class of $\Frob_2|\tK$ in $3.A_6$ is $(5ab)c^i$, where $c^i=(\P,\tK/K)^f= c^5=c^2$; so the coarse class of $\Frob_2|\tK$ is $15bd$.

\section{Summary of calculations}

We summarize here the calculations we have been able to complete. There are six cases in all for which we were able to test the conjecture. In the other cases, we determined that the level and weights were too big for computation at present.

\medskip

\renewcommand{\arraystretch}{1.1}
\begin{table}[hbt]
\begin{tabular}{l|l|l|l|l}
$A_6$ field   & $3.A_6$? & $p$  & $N,\eps$ &predicted weights \\  
\hline\hline
2,3: \#55     &yes       & 3    & $2^8$                  &(2,1,0)  \\ \hline
2,3: \#57     &yes       & 3    & $2^7$, $\psi_8$        &(5,3,1)   \\         \hline
2,3: \#58     &yes       & 3    & $2^7$, $\omega_4\psi_8$& (5,3,1) \\         \hline
3,7: \#3      &yes       & 3    & $7^2$    &(5,3,1)              \\ \hline
3,13: \#9     &no        & 3    & $13^2$   &(5,3,1)             \\ \hline
5,17: \#1     &yes       & 5    & $17$, $\eps_{17}$ &(6,3,0)   \\ \hline
\end{tabular}
\medskip
\caption{Completed cases}
\end{table} 

\medskip

In this table $\eps_{17}$ denotes the quadratic character with conductor $17$,  $\omega_4$ the quadratic character with conductor 4, and $\psi_8$ the real quadratic character of conductor 8. We have omitted the nebentype if it is 1.

Recall that when $p=5$ or $p=3$ there are 2 irreducible 3-dimensional representations of $3.A_6$ or $A_6$, respectively, each defined over $\FF_{p^2}$ and conjugate over $\FF_{p^2}$. Therefore, in each of the examples displayed in Table 4 the conjecture predicts 2 distinct Hecke eigenclasses, with eigenvalues in $\FF_{p^2}$ and conjugate over $\FF_{p^2}$,  to which the corresponding Galois representations appear to be attached, in the sense explained in the introduction.  This is indeed what we found. Of course, since the action of the  Hecke algebra on cohomology is defined over $\FF_p$, if we do find a Hecke eigenclass that appears to be attached to one of the two Galois representations, the conjugate eigenclass will appear to be attached to the other.

Note that verifying the $p=3$ examples requires determining which of the Frobenii of order 5 are conjugate in $A_6$.  We were able to do this in all cases except for $\Frob_2$ in the level $13^2$ example.  In that case all we could determine is that $\Frob_2$ has order 5. In this case, therefore,  we could not check the equality (2), though
our results do show that the left-hand and right-hand sides of (2) are either equal or conjugate over $\FF_3$. 

The first five of these examples take $p=3$, and so use a representation $\rho:G_\QQ \to GL_3(\bF_3)$ which cuts outs $K$; we don't need to know whether $K$ embeds in a $3.A_6$ extension of $\QQ$. But the case $3,13:\#9$, $p=3$ is interesting because there is in fact no $3.A_6$ extension: so the representation $\rho$  does not arise by reduction from a complex Galois representation. 

In the second through fifth  examples both the upper and the lower $2\times2$ block of $\rho|I_3$ are ``tr\`es ramifi\'ee,'' in the sense of section 2.4 (ii) of [Serre 87]. In accordance with the conjecture, we predict and find in each case an appropriate Hecke eigenclass in weight (5,3,1). But we have also checked that there are no such eigenclasses in weights $(3,3,1)$, $(3,1,1)$, or $(1,1,1)$, which are the other possible ways of resolving the ambiguity in the procedure used to determine the weights.

The sixth example has $p=5$, and therefore requires that $K$ can be extended to a $3.A_6$ extension of \QQ. This example was discussed in more detail in the previous section. 

\textbf{Remark:} The isomorphism $A_6\simeq PSL(2,\FF_9)$ gives rise  to a 2 dimensional linear representation $\sigma: 2.A_6\stackrel{\sim}{\rightarrow} SL(2,\FF_9)$ of the double cover $2.A_6$. The symmetric square  of this representation is trivial on the center of $2.A_6$, and in this way one can obtain the 3 dimensional representations of $A_6$ underlying the $p=3$ examples above. This raises the question of whether the representations $\rho:G_\QQ \to SL_3(\bF_9)$ in the first five examples in Table 4 are symmetric squares of representations $\tau: G_\QQ \to SL(2,\bF_9)$. If this were the case, then at least a weak form of our conjecture for these representations $\rho$ could be deduced from Serre's Conjecture ([Serre 87]; see also the discussion in [Ash-Sinnott 00], \S 3(i)).

However, this would require that the $A_6$ extensions in these examples embed in $2.A_6$ extensions. It is easy to see that they do not. All the twelve $A_6$ fields considered here are totally complex. If $K$ is a totally complex $A_6$ extension of $\QQ$, and $\hat K$ is a quadratic extension of $K$ such that $\Gal(\hat K/\QQ)\simeq 2.A_6$, let $c$ be a complex conjugation in $\Gal(\hat K/\QQ)$. Since the only element of $2.A_6$ of order 2 lies in the center, it follows that $c$ fixes $K$, which contradicts the fact that $K$ is totally complex.

\newpage

\hrulefill

\appendix
\section{Polynomials and extensions of $\QQ_p$}

The aim of this appendix is to find useful estimates for how accurately one needs to know the coefficients of an irreducible polynomial over a local field in order to have determined (up to isomorphism) the field extension obtained by adjoining a root. When factoring polynomials over $\QQ_p$, we need these estimates in order to determine the accuracy required in the factorization. 

The arguments are taken from the paper of Pauli and Roblot [Pauli-Roblot], adapted very slightly for our different purpose.

We work over $\QQ_p$, but the case of a general local field of characteristic 0 is no different. $|\cdot|$ denotes throughout the $p$-adic absolute value, normalized as usual by $|p|=1/p$.

Suppose that $f(x)=x^n+a_1 x^{n-1}+\cdots+a_n\in \ZZ_p[x]$ is an irreducible polynomial, and let $\alpha_1, \ldots,\alpha_n$ be the roots of $f$ (in $\bar\QQ_p$). Let 
\[
\delta_f=\min_{i\ne j} |\alpha_i-\alpha_j|
\]

We can calculate $\delta_f$ by finding the largest finite slope $\lambda$ in the Newton polygon of $f(x+\alpha_1)$. Then $\delta_f=|p|^\lambda$. A convenient way to calculate the Newton polygon of $f(x+\alpha_1)$ in PARI/GP is to work instead with $\tilde{f}(x)=\NN_{\QQ_p(\alpha_1)/\QQ_p}(f(x+\alpha_1))$; then $\tilde{f}$ lies in $\ZZ_p[x]$ and has the same set of slopes as does $f$ (the \emph{multiplicities} of these slopes, however, will have been multiplied by $n$).

We may also use the following lower bound for $\delta_f$:
\begin{lemma} Let $D_f=|disc(f)|$: then
    \[\delta_f\ge \left(\frac{D_f}{|a_n|^{n-2}}\right)^{1/n}.\]
\end{lemma}    
\begin{proof}
We suppose that the roots of $f$ have been numbered so that $\delta_f=|\alpha_1-\alpha_2|$. Then we have
\begin{align*}
|f'(\alpha_1)|&=\prod_{i=2}^n|\alpha_1-\alpha_i|\\
              &=\delta_f\prod_{i=3}^n|\alpha_1-\alpha_i|\\
	    &\le \delta_f\prod_{i=3}^n\max(|\alpha_1|,|\alpha_i|)\\
	    &\le \delta_f |a_n|^{(n-2)/n}
\end{align*}
Since $D_f=|f'(\alpha_1)|^n$, the above inequality may be written
\[
\delta_f\ge D_f^{1/n}/|a_n|^{(n-2)/n}
\]
\end{proof}

\bigskip

Let $f(x)=x^n+a_1 x^{n-1}+\cdots+a_n$ and $g(x)=x^n+b_1 x^{n-1}+\cdots+b_n$ be irreducible polynomials in $\ZZ_p[x]$; let $\alpha_1, \ldots,\alpha_n$ be the roots of $f$ and $\beta_1,\ldots, \beta_n$ be the roots of $g$. We consider the valuation of the resultant of $f$ and $g$:
\[
|R(f,g)|=\prod_{i,j=1}^n|\alpha_i-\beta_j|=\prod_{i=1}^n 
|g(\alpha_i)|=\prod_{j=1}^n|f(\beta_j)|.
\]
Since $f$ and $g$ are irreducible, this can be written
\[
|R(f,g)|=|g(\alpha_1)|^n=|f(\beta_1)|^n
\]

We have the following upper estimate for $|R(f,g)|$:

\begin{lemma}
    \[ |R(f,g)|^{1/n}\le|a_n|^{1/n}\max_{1\le i\le n} |b_i-a_i|\]
\end{lemma}
\begin{proof}
\begin{align*}
|R(f,g)|^{1/n}&=|g(\alpha_1)|\\
              &=|g(\alpha_1)-f(\alpha_1)|\\
              &=\left|\sum_{i=1}^n (b_i-a_i)\alpha_1^i\right|\\
              &\le\max_{1\le i\le n} |b_i-a_i||a_n|^{i/n}\\
              &\le|a_n|^{1/n}\max_{1\le i\le n} |b_i-a_i|\\
\end{align*}
(since $|\alpha_1|=|a_n|^{1/n}\le 1$).
\end{proof}

Next we derive an additional expression for $|R(f,g)|$. We suppose that the roots of $g$ have been numbered so that $|\alpha_1-\beta_1|$ is as small as possible. We have 
\[
|R(f,g)|=|g(\alpha_1)|^n=\prod_{j=1}^n|\alpha_1-\beta_j|^n.
\]
Now $|\alpha_1-\beta_j|\ge |\alpha_1-\beta_1|$ by our numbering; and we notice that 
\begin{itemize}
    \item if $|\alpha_1-\beta_j|>|\alpha_1-\beta_1|$, then 
\[
|\alpha_1-\beta_j|=|\alpha_1-\beta_j-(\alpha_1-\beta_1)|=|\beta_1-\beta_j|,
\]
so 
\[
|\alpha_1-\beta_j|=\max(|\alpha_1-\beta_1|, |\beta_1-\beta_j|);
\]
\item if  $|\alpha_1-\beta_j|=|\alpha_1-\beta_1|$, then 
\begin{align*}
|\beta_1-\beta_j|&=|\beta_1-\alpha_1+\alpha_1-\beta_j|\\
                 &\le \max(|\beta_1-\alpha_1||,|\alpha_1-\beta_j|)\\
	       &\quad=|\alpha_1-\beta_1|,
\end{align*}
so again we have
\[
|\alpha_1-\beta_j|=\max(|\alpha_1-\beta_1|, |\beta_1-\beta_j|);
\]
\end{itemize}
Thus we have
\begin{gather*}
|R(f,g)|^{1/n}=\prod_{j=1}^n \max(|\alpha_1-\beta_1|, |\beta_1-\beta_j|),
\intertext{ and in the same way, reversing the roles of $f$ and $g$: }
|R(f,g)|^{1/n}=\prod_{i=1}^n \max(|\alpha_1-\beta_1|, |\alpha_1-\alpha_i|),
\end{gather*}

We put these together to obtain a criterion for $f$ and $g$ to give the same extension of $\QQ_p$, i.e.\ for $\QQ_p[x]/(f(x))\simeq \QQ_p[x]/(g(x))$:

\begin{prop}   If $f$ and $g$ are irreducible and if
\begin{gather*}
\max_{1\le i\le n} |b_i-a_i|<\delta_f \left(\frac{D_f}{|a_n|}\right)^{1/n} \\
\intertext{ or if }
\max_{1\le i\le n} |b_i-a_i|< 
\left(\frac{D_f^2}{|a_n|^{n-1}}\right)^{1/n} 
\end{gather*}

---then $f$ and $g$ give the same extension of $\QQ_p$.
\end{prop}
\begin{proof}

Assume that $f$ and $g$ give different extensions of $\QQ_p$; then by Krasner's Lemma, \[
\delta_f\le|\alpha_i-\beta_j|.
\]
for all $i$ and $j$.

We give a lower bound on $|R(f,g)|$ based on this inequality.

\begin{align*}
 |R(f,g)|^{1/n}&=\prod_{i=1}^n \max(|\alpha_1-\beta_1|, |\alpha_1-\alpha_i|)\\
                &\ge \prod_{i=1}^n \max(\delta_f, |\alpha_1-\alpha_i|)\\ 
                &\quad=\delta_f\prod_{i=2}^n \max(\delta_f, |\alpha_1-\alpha_i|)\\ 
                &\quad=\delta_f\prod_{i=2}^n|\alpha_1-\alpha_i|\\
                &\quad=\delta_f|f'(\alpha_1)|\\
                &\quad=\delta_f D_f^{1/n}.\\
\end{align*}  
So the first inequality in the statement of the proposition leads to a contradiction with Lemma 2. The second inequality implies the first, by Lemma 1, and so leads to a contradiction as well.
\end{proof}

\newtheorem*{lembis}{Lemma 1 bis}
\newtheorem*{propbis}{Proposition 1 bis}

Finally, we rewrite these results using the $p$-adic valuation $\ord_p$. Let $d=\ord_p disc(f)$; let $a=\ord_p(a_n)$; and let $\lambda=\max_{i\ne j} \ord_p(\alpha_i-\alpha_j)$. Then we have

\begin{lembis}\[\lambda\le \frac{d-(n-2)a}{n}.\]
\end{lembis}

\begin{propbis} If $k$ is an integer such that
\[
k>\lambda+\frac{d-a}{n}
\]
or such that
\[
k>\frac{2d-(n-1)a}{n},
\]
Then any monic irreducible polynomial of degree $n$ congruent to $f \bmod p^k$ gives the same extension of $\QQ_p$ as $f$. 
\end{propbis}

\emph{Caveat:} The second polynomial \emph{must} be irreducible for the result to apply. If we simply calculate $k$ and find a monic polynomial $g\in\ZZ_p[x]$ such that $f\equiv g \bmod p^k$, it can happen that $g$ is reducible and fails to give the same extension of $\QQ_p$. For example, let $p=2$ and let $f=x^3-2$; then $d=2$, $a=1$, so that $(2d-(n-1)a)/n=2/3$ and we may take $k=1$. We cannot take $g=x^3$, but the Proposition does say that any irreducible polynomial $g\equiv x^3 \bmod 2$ gives the same extension as $f$.

\medskip

Avner Ash,
Department of Mathematics,
Boston College\\
\texttt{ashav@bc.edu}

David Pollack,
Department of Mathematics and Computer Science,
Wesleyan University\\
\texttt{dpollack@wesleyan.edu}

Warren M.\ Sinnott,
Department of Mathematics,
Ohio State University\\
\texttt{sinnott@math.ohio-state.edu}

\end{document}